\documentclass[11pt,reqno]{amsart}
\usepackage{fullpage,url,amssymb,enumerate,colonequals}\usepackage{float}
\usepackage[usenames]{color}
\usepackage{amssymb}
\usepackage{amsmath}
\usepackage{amsthm}
\usepackage{amsfonts}
\usepackage{amscd}
\usepackage{graphicx}
\usepackage{enumitem}
\usepackage{changepage}
\usepackage{caption}
\usepackage{tikz}
\usetikzlibrary{trees}
\usepackage{verbatim}

\usepackage{hyperref}
\usepackage{color}
\usepackage{fullpage}
\usepackage{float}
\usepackage{graphics}
\usepackage{latexsym}
\usepackage{epsf}
\usepackage{breakurl}
\usepackage[english]{babel}
\DeclareMathOperator{\lcm}{lcm}

\usepackage[left=1in,right=1in,top=1in,bottom=1in]{geometry}

\theoremstyle{plain}
\newtheorem{theorem}{Theorem}
\newtheorem{corollary}[theorem]{Corollary}
\newtheorem{lemma}[theorem]{Lemma}
\newtheorem{proposition}[theorem]{Proposition}

\theoremstyle{definition}

\newtheorem{example}[theorem]{Example}

\theoremstyle{remark}

\newcommand{\seqnum}[1]{\href{https://oeis.org/#1}{\rm \underline{#1}}}
\makeatletter
\@namedef{subjclassname@2020}{%
  \textup{2020} Mathematics Subject Classification}
\makeatother

\begin{document}
\title{Linear Recurrences of Order at most two in Small Divisors}

\author{A. Anas Chentouf}
\address{Department of Mathematics, Massachusetts Institute of Technology, Cambridge, MA 02139, USA}
\email{chentouf@mit.edu}
\begin{abstract}
\large Given a positive integer $n$, the small divisors of $n$ are defined as the positive divisors that do not exceed $\sqrt{n}.$ Ianucci  previously classified all $n$ for which  the small divisors of $n$ form an arithmetic progression. In this paper, we classify all $n$ for which the small divisors of $n$ form a linear recurrence of order at most two. 
\end{abstract}

\maketitle

\section{Introduction}
\large
\par
As usual, we say that a non-zero integer $m$ divides another integer $n$ whenever the quotient $\frac{n}{m}$ is itself an integer. When this occurs, we use the conventional notation $m|n$. Moreover, a \textit{nontrivial} divisor of a natural number $n$ is a divisor other than $1$ and $n$ (if such a divisor exists).
\par
Consider a positive integer $n,$ and let $\mathcal{S}_n$ be the set of \textit{small} divisors of $n,$ that is, 
\begin{equation}
    \mathcal{S}_n:= \{ d: 1 \le d \le \sqrt{n} \land d|n \}. \label{smallset}
\end{equation} \textit{Large} divisors, associated with $\mathcal{L}_n,$ are defined with the inequality reversed. For brevity, we let $s(n)=|\mathcal{S}_n|$.  Moreover, for any tuple $(u,v,a,b) \in \mathbb{Z}^4,$ there is an integral linear recurrence of order at most two, $(d_i)_{i=1}^\infty,$ given by 

\begin{equation}
d_i = \begin{cases}
 u, & \text{ if } i=1; \\
v, & \text{ if } i=2; \\
ad_{i-1}+bd_{i-2}, & \text{ if } \ i \ge 3.
\end{cases} \label{recurrences}
\end{equation}

This recurrence is usually abbreviated as $U(u,v,a,b).$ The choice of beginning the index at $1,$ rather than $0,$ is due to the context of working with divisors.

As usual, we let $\tau(n)$ denote the number of positive divisors of $n,$ that is, 

$$
    \tau(n) = \sum_{d | n}^{} 1.
$$

Since there is a ‘‘coupling’’ between small and large divisors of $n,$ linking $d$ and $\frac{n}{d},$ we easily relate $\tau(n)$ to $s(n)$ as follows:

\begin{equation}
\label{tau}
\tau(n) =
    \begin{cases}
    
        2 s(n)-1, & \text{if $n$ is a perfect square;} \\
         2 s(n),& \text{otherwise.}
    \end{cases}
\end{equation}
In this paper, we classify all values of $n \in \mathbb{N}$ such that the \textit{ordered} elements of $\mathcal{S}_n$ form an integral linear recurrence $ (d_i)_{n=1}^{s(n)}.$ A number $n$ that satisfies this property is called a \textit{recurrent} number. 

\begin{example}
Notice that $60$ is a recurrent number since the elements of  $\mathcal{S}_{60}=\{1,2,3,4,5,6\}$ satisfy the recurrence $d_{i}=2d_{i-1}-d_{i-2}.$ On the other hand, a simple computation or computer search would show that the smallest non-recurrent number is $36$,  as  $\mathcal{S}_{36}=\{1,2,3,4,6\}$ does not satisfy a  linear recurrence of order at most two.
\end{example}

Since $1$ is the smallest divisor of $n,$ we have $d_1=1.$ To be able to easily determine whether $\mathcal{S}_n$ forms a linear recurrence, we  impose the condition that $s(n) \ge 5,$ although we revisit the cases $s(n) \le 4$ at the end of the paper. Furthermore, since the smallest nontrivial divisor of any natural number is prime, we henceforth assign $d_2=p$ throughout the paper. We also let $U(1,p,a,b)$ be the linear recurrence which produces the divisors of some recurrent number $n.$ Notice that simply knowing $p,a,b$ allows us to determine all element in $\mathcal{S}_n$ whenever $n$ is recurrent.

\begin{center}
\begin{equation}
  \begin{tabular}{|c|c|}
  \hline
     $d_1$ & $1$ \\  \hline
     $d_2$ & $p$ \\  \hline
     $d_3$ & $ap+b$ \\ \hline
     $d_4$ & $(a^2+b)p+ab$ \\ \hline
     $d_5$ & $(a^3+2ab)p+b(a^2+b)$ \\ \hline
\end{tabular}  
\end{equation}
\captionof{table}{First $5$ values of $(d_i)_{i=1}^{s(n)}$ in terms of $a,b,p.$}
\label{sophisticatedtable}
\end{center}

\tikzstyle{level 1}=[level distance=3.5cm, sibling distance=4cm]
\tikzstyle{level 2}=[level distance=2cm, sibling distance=4cm]
\tikzstyle{level 3}=[level distance=3.5cm, sibling distance=2cm]

\tikzstyle{bag} = [text width=2em, text centered]
\tikzstyle{end} = [circle, minimum width=2pt,fill, inner sep=1pt]

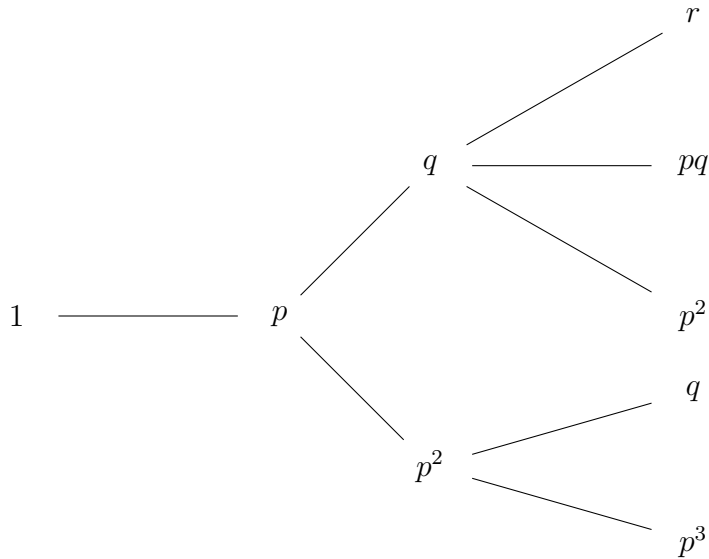
\begin{figure}[H]
    \centering

\begin{tikzpicture}[grow=right, sloped]
\node[bag] {1}
    child {
        node[bag] {$p$}        
            child {
                node[bag]{$p^2$} {}
                child {
                    node[bag]{$p^3$} {}
                    edge from parent
                    node[above] {}
                    node[below]  {}
                }
                child {
                    node[bag]{$q$} {}
                    edge from parent
                    node[above] {}
                    node[below]  {}
                }
                edge from parent
                node[above] {}
                node[below]  {}
            }
            child {
                node[bag]{$q$} {}
                child {
                    node[bag]{$p^2$} {}
                    edge from parent
                    node[above] {}
                    node[below]  {}
                }
                child {
                    node[bag]{$pq$} {}
                    edge from parent
                    node[above] {}
                    node[below]  {}
                }
                child {
                    node[bag]{$r$} {}
                    edge from parent
                    node[above] {}
                    node[below]  {}
                }
                edge from parent
                node[above] {}
                node[below]  {}
            }
            edge from parent 
            node[above] {}
            node[below]  {}
    };
    
\end{tikzpicture}
    \caption{The above tree represents the possible ‘‘configurations’’ of the first four divisors of a natural number $n \ge 1,$ where $p<q<r$ are distinct primes. For example, this says that the first divisor must be $1$, the second divisor must be some prime $p$, and the third divisor can either be $p^2$ or a second prime $q$, etc... Our approach makes use of this "tree" approach in determining the divisors $d_i.$}
    \label{fig:my_label}
\end{figure}

The results obtained by Ianucci \cite{Ianucci} appear as special cases in this generalized characterization. This is because any arithmetic progression satisfies the second-order linear recurrence $d_{i}=2d_{i-1}-d_{i-2}.$ Many recurrent numbers will have small divisors forming a geometric sequence $d_{i}=ad_{i-1},$ or a \textit{bifurcated} geometric sequence, $d_{i}=bd_{i-2}.$  However, unlike with arithmetic progressions, we will see that for recurrent numbers $n,$ the value of $s(n)$ may be arbitrarily large. This changes the nature of the approach adopted in this paper. We begin by reducing cases on $a,b$ and working with the values of $\gcd(a,b),$ as well as using the aforementioned tree approach. 
\section{Preliminary Results}

\begin{proposition}
\label{div}
Let $a,b,n$ be positive integers such that $a \le b$. If $ab|n,$ then $a \in \mathcal{S}_n.$ 
\end{proposition}
\begin{proof}
Since $ab|n,$ it follows that $a$ is a divisor of $n$. Moreover, we get that $ab \le n.$ However, $a^2 \le ab \le n,$ and hence $a \le \sqrt{n}.$
\end{proof}

\
\begin{theorem}
\label{thm}
Given a recurrent number $n$ associated with $U(1,p,a,b),$ at least one of the following statements is true:
\begin{enumerate}
    \item $\gcd(a,b)=1.$
    \item All nontrivial elements of $\mathcal{S}_n$ are divisible by $p.$
    \item $a=0.$
\end{enumerate}
\end{theorem}
\begin{proof}
Suppose $\gcd(a,b) \ne 1,$ so some prime divides both $a$ and $b$. 

If $p|\gcd(a,b),$ then we see that all nontrivial elements in $\mathcal{S}_n$ are divisible by $p,$ so the second possibility holds. If $p \nmid \gcd{(a,b)},$ then we obtain that $\gcd{(a,b)}|ap+b=d_3$ so $\gcd{(a,b)}$ is some prime $q,$ and $d_3=q.$ Similarly, $q|d_4$ and thus $aq+bp=d_4$ must be $pq.$ Equating, we get $a(p^2-q)=0,$ hence $a=0.$
\end{proof}

\begin{proposition}
\label{pkq}
If all nontrivial elements of $\mathcal{S}_n$ are divisible by some prime $p,$ then $n=p^k$ or $n=p^kq$ for some prime $q>p^k.$ 
\end{proposition}
\begin{proof}
We proceed by considering two cases.

\begin{adjustwidth}{0.5cm}{}
Case I: If all nontrivial divisors of $n$ are divisible by $p,$ then $n=p^k$ for some $k \ge 1.$
\end{adjustwidth}
\begin{adjustwidth}{0.5cm}{}
Case II: Assume that $n$ has a divisor $q> \sqrt{n}$ that is not divisible by $p.$ 
\begin{adjustwidth}{0.5cm}{}

If $q$ is composite, then one of its prime divisors must be less than or equal to $\sqrt{q}< \sqrt{n},$ and must hence be divisible by $p.$ This contradicts our assumption, and $q$ must thus be prime. 
\end{adjustwidth}
\begin{adjustwidth}{0.5cm}{}

We now prove uniqueness. If there exists another prime $r$ in the interval $(\sqrt{n}, n)$ that divides $n$, then $qr$ must divide $n$ but $qr > \sqrt{n} \cdot \sqrt{n} =n,$ another contradiction.

Hence, $q$ is the sole divisor of $n$ that is not divisible by $p,$ and thus we have that $n=p^kq$.
\end{adjustwidth}

\end{adjustwidth}
 
\end{proof}
\begin{proposition}
In a recurrent number $n$, if $\mathcal{S}_n$ satisfies a recurrence $U(1,p,0,b)$ (i.e., $a=0$), then exactly one of the following possibilities holds, where  $p,q,r$ are primes and $k \in \mathbb{N}$:

\begin{enumerate}
    \item $n=p^k$ 
    \item $n=pq^k$ for $p<q.$
    \item $n=pq^k r$ with $p<q<pq^{k}<r.$ 
\end{enumerate}

\end{proposition}
\begin{proof}
If $a=0,$ then $\mathcal{S}_n$ forms a bifurcated series with $d_i=bd_{i-2},$ but taking $i=3$ requires that $d_3=b$ must be either $p^2$ or $q$ (for a prime $q>p$). The first case gives rise to $\mathcal{S}_n=\{1,p,p^2, \dots , p^{r}\}$, while the second implies that  $\mathcal{S}_n=\{1,p,q,pq,q^2 \dots\ pq^l\}$ forms a bifurcated sequence. Identical to the proof of Lemma \ref{pkq}, there is at most one prime divisor greater than $\sqrt{n}$ that divides $n,$ and hence the possibilities follow directly.
\end{proof}

\section{The case $\gcd(a,b)=1$}

We now turn our attention to the major case where $\gcd(a,b)=1.$ We also assume, henceforth, that $n$ has at least two distinct small prime factors - otherwise we are simply in the second and third cases of Theorem \ref{thm}. In particular, $ab \ne 0$ since if one of the numbers is zero, then the other cannot have absolute value $1$, as the sequence of small divisors cannot negative numbers and cannot be constant. 

\begin{lemma}
\label{bco}
For all small divisors $d_i$ of a recurrent number, $\gcd({b,d_i})=1.$ 
\end{lemma}
\begin{proof}
Suppose that  $q|\gcd({b,d_i})$ for some prime $q$ and some index $i.$ Notice that this $i \ne 1$. Similarly, if $i=2,$ then $p|b$ and inductively, all small divisors are divisible by $p$ - this contradicts the assumption of having two distinct small prime factors.

Otherwise, $i \ge 3$ and $d_i=ad_{i-1}+bd_{i-2},$  hence $q|ad_{i-1}.$ The coprimality of $a$ and $b$ implies that $q|d_{i-1}.$ We can inductively show that smaller divisors, including $d_1=1$, are divisible by $q,$ which is absurd. Hence the coprimality follows.
\end{proof}

\begin{corollary}
\label{coprime}
For all $i$ such that $d_i, d_{i+1}$ are small divisors of a recurrent number, we have $\gcd({d_i, d_{i+1}})=1.$
\end{corollary}
\begin{proof}
We proceed by induction. The base case clearly holds. For the inductive step, observe that
$$\gcd{(d_i, d_{i+1}})=\gcd{(d_i, ad_{i}+bd_{i-1})} = \gcd{(d_i, bd_{i-1})}.$$

The inductive hypothesis implies that $\gcd{(d_i, d_{i-1})}=1,$ while Lemma \ref{bco} implies that $\gcd{(d_i,b})=1,$ and hence the coprimality follows.
\end{proof}

As a direct consequence of this, we conclude that $d_3$ cannot be $p^2$ in a recurrent number, and hence $d_3$ must be prime. Henceforth, we set $d_3=q.$

\begin{lemma}
All configurations of the first four non-trivial divisors $(d_2,d_3,d_4,d_5)$ except possibly $(p,q,p^2,r), (p,q,r,p^2)$ and $(p,q,r,s)$, where $p,q,r,s$ are distinct primes, are impossible in a recurrent number.
\end{lemma}
\begin{proof}
We have already determined $d_2,d_3.$ Now, $d_4$ can be $p^2,q^2,pq,$ or $r$ (where $r$ is another prime). By Corollary \ref{coprime}, it cannot be $pq$ or $q^2$. Hence, we have either $d_4=p^2$ or $d_4=r.$

If $d_4=p^2,$ then we can similarly see that the possibilities for the next divisor $d_5$ are $p^3,q^2,pq,r$ for a prime $r.$ However, Corollary \ref{coprime} once again implies that $p^3, pq$ are impossible. Similarly speaking, if $d_5=q^2,$ then since $d_5=ap^2+bq$ we would get that $q|ap^2,$ and so, $q|a.$ Yet, since $q=ap+b, q|a$ would imply that $q|b$ and thus contradicts Lemma \ref{coprime}.  Hence, $d_5$ ought to be another prime, $r.$

If $d_4=r,$ then we see that our possibilities for the next divisor $d_5$ are $p^2,pq,s,$ for a new prime $s.$  A reasoning similar to above shows that $q \nmid d_5,$ so the only possibilities left are those of $d_5=p^2$ and $d_5=s.$ 
\end{proof}

\begin{lemma}
\label{bounds}
The only pair $(u,x) \in \mathbb{N} \times \mathbb{N}_{\ge 2}$ that satisfies 
$$-ux^5+(u+3)x^4-\frac{2u+3}{u}x^3+\frac{3u+1}{u^2}x^2-\frac{2}{u^2}x+\frac{1}{u^2} \ge x^2+1$$
is $(1,2).$
\end{lemma}
\begin{proof}
Define $$P(x):=-ux^5+(u+3)x^4-\frac{2u+3}{u}x^3+\frac{3u+1}{u^2}x^2-\frac{2}{u^2}x+\frac{1}{u^2} - (x^2+1).$$ Notice that 
\begin{align*}
    P'(x) &=-5ux^4+4(u+3)x^3-\frac{3(2u+3)}{u}x^2+\frac{2(3u+1)}{u^2}x-\frac{2}{u^2}-2x \\
    &= -5ux^4+4(u+3)x^3-\frac{3(2u+3)}{u}x^2+\frac{-2u^2+6u+2}{u^2}x-\frac{2}{u^2}
\end{align*}
$$$$

For $u \ge 4, x \ge 2$ we have that 
$$5ux^4 \ge 10ux^3 \ge 4(u+3)x^3,$$ and in addition, we  have $-\frac{3(2u+3)}{u}x^2,\frac{-2u^2+6u+2}{u^2}x,-\frac{2}{u^2}$ are each negative. Hence, $P'(x)$ is negative and for $u \ge 4,$ the function $P(x)$ is strictly decreasing in the interval $[2, \infty).$ Moreover, it is easy to check that $P(2)<0.$ Hence, it suffices to simply consider $u \in \{1,2,3\}.$ One may easily check that the only possible value is $u=1$. Another simple computation reveals that the inequality only holds for $x=2.$
\end{proof}
\begin{theorem}
The only recurrent number with configuration  $(d_2,d_3,d_4,d_5)=(p,q,p^2,r)$ is $60.$
\end{theorem}
\begin{proof}
In this case, we have that 

\begin{equation}
\label{kp}
p^2=d_4=aq+bp=a(ap+b)+bp.
\end{equation}

Since $p|ab,$ yet by Lemma \ref{bco}, $p \nmid b$ and hence $p|a.$ Hence,  $a=kp$ for some integer $k.$ We can rewrite (\ref{kp}) as 

$$k^2p^2 +b(k+1)=p,$$

which is equivalent to 
$$
b = \frac{p(1-k^2p)}{k+1}.$$ However, since we know that $p \nmid b,$ we obtain that $p|k+1.$ In other words, there exists an integer $u$ such that $k+1=up.$ From here, we are able to express $a,b$ in terms of $u,p,$ and use that to obtain $r$ in terms of $u,p$ only. We end up obtaining

$$a=up^2-p,$$
$$b=\frac{1-p}{u}-up^3+2p^2,$$
and
$$r=ap^2+bq=-up^5+(u+3)p^4-\frac{2u+3}{u}p^3+\frac{3u+1}{u^2}p^2-\frac{2}{u^2}p+\frac{1}{u^2}.$$
Similarly, we note that 
$$u=\frac{p-1}{p^2-q} > 0.$$ Since $r=d_5\ge d_4+1 = p^2+1$, Lemma \ref{bounds} gives that $p=2,$ and substituting reveals $a=2, b=-1,$  and hence the first five divisors would be $1,2,3,4,5.$ This  reduces to the case of an arithmetic sequence with common difference $1,$ which, by Lemma 2 of \cite{Ianucci}, implies that $s(n)=5$ or $s(n)=6.$ Having 5 divisors would still require that $n \ge \lcm[1,2,3,4,5]=60,$ but $6$ is a small divisor that does not appear in this set, but this contradicts the recurrency of $n.$ When $s(n)=6,$ we have that $n$ is a multiple of $60.$ However, if $n \ge 180$ then $12 \le \sqrt{n}$ does not appear in the list of small divisors, contradicting the recurrency. Clearly, $n=120$ is not recurrent either as $8 \notin \{1, 2, 3,4,5,6\}$. Hence, the only solution is  $\mathcal{S}_n=\{1,2,3,4,5,6\}$ for when $n=60,$ which is indeed recurrent. 
\end{proof}

\begin{corollary}
\label{three}
For an initial configuration $(p,q,r,p^2)$ in a recurrent number, we have that $p|d_{j}$ if and only if $j \equiv 2 (\bmod{3}).$
\end{corollary}

\begin{proof}
Firstly, the statement holds for $d_2=p. $ Since $p^2=d_5=(a^3+2ab)p+b(a^2+b),$ Lemma \ref{bco} implies that $p|(a^2+b).$ Hence, $$d_{i+3}=ad_{i+2}+bd_{i+1}=a(ad_{i+1}+bd_i)+bd_{i+1}=(a^2+b)d_{i+1}+abd_{i}.$$ Lemma \ref{bco} once again gives us that $p|d_{i+3}$ if and only if $p|d_{i},$ which is the inductive step.
\end{proof}

\begin{theorem}
\label{impossible}
The configuration of divisors $(p,q,r,p^2)$ cannot occur in a recurrent number.
\end{theorem}
\begin{proof}
Assume otherwise, that is, there exists a recurrent number $n$ of this configuration. We prove, by induction, that $s(n) \ge 3i+2$ for all $i \in \mathbb{N}.$ This would be absurd as $s(n)$ grows endlessly.

The base case is already assumed. Assume that $s(n) \ge 3i+2$ for some $i \in \mathbb{N}.$ By Corollary \ref{three}, $p\nmid d_{3i}d_{3i+1}$ and since $\gcd(d_{3i},d_{3i+1})=1$ by Corollary \ref{coprime}, we have that $p^2d_{3i}d_{3i+1}$ divides $n.$ By Proposition \ref{div}, $pd_{3i}$ is a small divisor. We consider two cases.

\begin{adjustwidth}{0.5cm}{}
Case I: $pd_{3i}=d_{3i+2}.$ Hence, we obtain that $d_{3i}|ad_{3i+1},$ and by Lemma $\ref{coprime}$, we obtain that $d_{3i}|a.$ Observe that the sequence $(d_i \bmod{a})$ takes values congruent to $1,p,b,bp,b^2, \dots,$ and thus, $d_{3i}$ divides an element of the form $b^kp$ where $k \in \mathbb{N}.$ However, the assumption and Lemma \ref{bco}  yield $$\gcd(p,d_{3i+1})=1=\gcd(b,d_{3i+1}),$$ which is a contradiction.  
\end{adjustwidth}
\begin{adjustwidth}{0.5cm}{}
Case II: $pd_{3i}>d_{3i+2},$ hence, $pd_{3i} \ge d_{3i+5}$ by Corollary \ref{three}. Therefore, $s(n) \ge 3i+5,$ and the induction is complete.  
\end{adjustwidth}

\end{proof}

\begin{lemma}
\label{lastlem}
In a recurrent number of configuration $(p,q,r,s)$, the small divisors $d_i,d_{i+2}$ are coprime.
\end{lemma}
\begin{proof}
Assume a prime $t$ divides $\gcd{(d_{i}, d_{i+2})}.$ Thus, we have that $t|ad_{i+1}.$ Clearly, $t \nmid d_{i+1}$ by Lemma \ref{coprime}, and therefore, $t|a.$ As such, we have that the sequence $(d_i \bmod{t})$ takes values congruent to $1,p,b,bp,b^2,\dots$ and hence we have that $t|b^kp$ for some $k \ge 0.$ Yet, Lemma \ref{bco} implies that $\gcd{(b,t)}=1,$ hence $t|p$ and so $t=p.$ However, this would imply that $t|r=(a^2+b)p+ab,$ which is a contradiction as $p<r$ are distinct primes.
\end{proof}

\begin{theorem}
The configuration of divisors $(p,q,r,s)$ cannot occur in a recurrent number.
\end{theorem}
\begin{proof}
Assume that it does occur in some recurrent number $n$. Since $pq<pr$ are both small by Proposition (\ref{div}), there exist at least three small divisors that are multiples of $p.$ Let $d_{j}=pv$ be the greatest small divisor that is divisible by $p$. In particular, $v > p.$ By Lemma \ref{lastlem}, $p \nmid d_{j-2}d_{j-1},$ and we conclude that $d_{j-2}d_{j-1}pv|n,$ so Proposition \ref{div} implies that $pd_{j-2}\in \mathcal{S}_n$. Moreover, $pd_{j-2} \ne d_{j}$ by Lemma \ref{lastlem}, and so $pd_{j-2}>d_{j}$ is a small divisor greater than $d_j,$ a contradiction. 

\end{proof}

\section{Concluding Remarks}

In conclusion, imposing the condition $s(n) \ge 5$ returns certain infinite families of recurrent integers $n$ with at most $3$ distinct prime divisors, in addition to the sole ‘‘sporadic’’ case of $60.$ We now revisit the condition we imposed and relax it to find all recurrent numbers by relating $s(n)$ to $\tau(n)$ using Equation (\ref{tau}).

\begin{enumerate}[label=(\roman*)]
\item $s(n)=1$. Hence, $\tau(n) \in \{1,2\}.$ That is, $n=1$ or $n=p$ for some prime $p,$ and in both cases, $n$ is vacuously recurrent.
\item $s(n)=2$. Thus, $\tau(n) \in \{3,4\},$ and we have the possibilities that $n=p^2,p^3$, or $pq$ for primes $p<q.$ In all cases, $n$ is vacuously recurrent.
\item $s(n)=3$. This implies $\tau(n) \in \{5,6\}.$ and hence the possibilities are $n=p^4$, $n=p^5$, $n=p^2q$ or $n=pq^2$ for primes $p<q.$ These are again vacuously recurrent. 
\item $s(n)=4$. Consequently, $\tau(n) \in \{7,8\}.$ 

\begin{enumerate}[label=(\Roman*)]
    \item If $\tau(n)=7,$ we get that $n=p^6$ and $\mathcal{S}_n=\{1,p,p^2,p^3\}$ indeed forms a linear recurrence of order at most two. 
    \item If $\tau(n)=8,$ the possibilities are $p^7,p^3q,$, $pq^3$ and $pqr$ for primes $p<q<r.$ 
    \begin{enumerate}
        \item $n=p^7$: Here, $n$ is recurrent by the aforementioned geometric recurrences. 
        \item $n=p^3q$: Here, Lemma \ref{coprime} implies that the only possibility is $\mathcal{S}_n= \{1,p,q,p^2\}$. Setting up a linear system of equations to solve for $a,b$, we deduce that it is both necessary and sufficient for $\frac{q-p}{p^2-q}$ to be an integer for primes $p<q<p^2$. Based on numerical evidence, we conjecture that this occurs for infinitely many pairs of primes $(p,q)$, but we defer analyzing this conjecture to another time.
        \item $n=pq^3:$ Clearly, we have that $\mathcal{S}_n=\{1,p,q,pq\},$ but this is not possible by Lemma \ref{coprime}. 
        \item $n=pqr:$ If $pq<r,$ then $\mathcal{S}_n=\{1,p,q,pq\}$ and the same argument above holds. If $pq>r,$ then $\mathcal{S}_n=\{1,p,q,r\}.$ Setting up a linear system of equations, we see that it is both necessary and sufficient for $\frac{pq-r}{p^2-q}$ to be an integer.
    \end{enumerate}
    
\end{enumerate}

\end{enumerate}

We summarize our findings in the following result.

\begin{theorem}
All recurrent numbers $n$ fall into one of the following categories.

\begin{enumerate}
\item $n=p^k:$ for some prime $p$ and a natural number $k,$ hence $\mathcal{S}_n=\{1,p, \dots ,p^{\lfloor \frac{k}{2}\rfloor}\}$. 
\item $n=p^kq:$ for some primes $p,q$ and a natural number $k,$ such that $q>p^k$ . Hence, $\mathcal{S}_n=\{1,p, \dots ,p^{k}\}.$ 
\item $n=pq^k:$ for some primes $p<q,$  and $k$ an odd (resp.,  even) natural number with $ \mathcal{S}_n=\{ 1,p,q, pq, q^2, \dots ,pq^{\frac{k-1}{2}} \}$ (resp.,  $\mathcal{S}_n=\{ 1,p,q, pq, q^2, \dots ,q^{\frac{k}{2}} \}$). 
\item $n=pq^kr:$ for some primes $p,q,r$ and a natural number $k$ such that $p<q$ and $r>pq^k,$ hence, $\mathcal{S}_n=\{1,p,q,pq, \dots ,pq^k\}.$
\item $n=60:$ with $\mathcal{S}_n=\{1,2,3,4,5,6\}.$
\item $n=p^3q$ for some primes $p,q$ such that $p<q<p^2$ and $p^2-q|q-p.$ Hence, $\mathcal{S}_n=\{1,p,q,p^2\}.$
\item $n=pqr$ for some primes $p,q,r$ such that $p<q<r$ and $p^2-q|pq-r.$ Hence, $\mathcal{S}_n=\{1,p,q,r\}.$
\end{enumerate}
\end{theorem}

Notice that $\mathcal{S}_n$ may contain an arbitrarily long subset of divisors that form a second-order linear recurrence, without all of $\mathcal{S}_n$ forming such a recurrence. Consider numbers of the form $pq^kr^2$ where $p<q<r$ are primes and $k \in \mathbb{N}$ such that $pq^k<r.$ If $r$ is large enough, we see that the divisors $1,p,q,pq,q^2, \dots $ may form an arbitrarily long linear recurrence that is interrupted by $r,$ but $n$ is not itself recurrent.

Finally, we briefly study the analytic distribution of recurrent numbers. Define $f(x)$ to be the number of  recurrent integers in the interval $[1,x].$ The fact that recurrent numbers have at most $3$ distinct prime divisors implies that they have null density among the natural numbers. For $k \in \mathbb{N},$ let $\pi_k(x)$ be the number of integers in $[1,x]$ having exactly $k$ distinct prime factors. Hardy and Ramanujan \cite{Ramanujan} showed that there exist constants $A,B$ such that
$$\pi_k(x) < \frac{Ax(\log\log{x}+B)^{k-1}}{(k-1)!\log{x}}.$$ In fact, Landau \cite{Landau} proved that that $$\pi_k(x) \sim \frac{x(\log\log{x})^{k-1}}{(k-1)!\log{x}}.$$

Since $f(x) \le \pi_1(x)+\pi_2(x)+\pi_3(x)$, we obtain that there exists some constant $C$ such that for all $x \ge 1,$
$$f(x)\le C \Bigg(\frac{x}{\log{x}}+\frac{x\log{\log{x}}}{\log{x}}+\frac{x(\log{\log{x}})^2}{2\log{x}}\Bigg).$$

\section{Acknowledgment}
I would like to thank Professor Haynes Miller for his feedback and support, and Murilo Curato Zanarella for guidance in dealing with the analytic distributions. I would also like to thank Abdeljalil Hezouat for his comments.

\bigskip
\hrule
\bigskip

\noindent
(Concerned with OEIS sequence
\seqnum{A346447}.)

\bigskip
\hrule
\bigskip

\vspace*{+.1in}
\noindent

\end{document}